# PREDICTION IN FUNCTIONAL LINEAR REGRESSION

### By T. Tony Cai[1] and Peter Hall

*University of Pennsylvania and Australian National University*

There has been substantial recent work on methods for estimating the slope function in linear regression for functional data analysis. However, as in the case of more conventional finite-dimensional regression, much of the practical interest in the slope centers on its application for the purpose of prediction, rather than on its significance in its own right. We show that the problems of slope-function estimation, and of prediction from an estimator of the slope function, have very different characteristics. While the former is intrinsically nonparametric, the latter can be either nonparametric or semiparametric. In particular, the optimal mean-square convergence rate of predictors is $n^{-1}$, where $n$ denotes sample size, if the predictand is a sufficiently smooth function. In other cases, convergence occurs at a polynomial rate that is strictly slower than $n^{-1}$. At the boundary between these two regimes, the mean-square convergence rate is less than $n^{-1}$ by only a logarithmic factor. More generally, the rate of convergence of the predicted value of the mean response in the regression model, given a particular value of the explanatory variable, is determined by a subtle interaction among the smoothness of the predictand, of the slope function in the model, and of the autocovariance function for the distribution of explanatory variables.

**1. Introduction.** In the problem of functional linear regression we observe data $\{(X_1, Y_1), \ldots, (X_n, Y_n)\}$, where the $X_i$'s are independent and identically distributed as a random function $X$, defined on an interval $\mathcal{I}$, and the $Y_i$'s are generated by the regression model,

$$(1.1) \qquad Y_i = a + \int_{\mathcal{I}} b\, X_i + \varepsilon_i.$$

Received August 2004; revised October 2005.

[1]Supported in part by NSF Grant DMS-03-06576 and a grant from the Australian Research Council.

Here, $a$ is a constant, denoting the intercept in the model, and $b$ is a square-integrable function on $\mathcal{I}$, representing the slope function. The majority of attention usually focuses on estimating $b$, typically by methods based on functional principal components. See, for example, [28], Chapter 10, and [29].

In functional linear regression, perhaps as distinct from more conventional linear regression, there is significant interest in $b$ in its own right. In particular, since $b$ is a function rather than a scalar, then knowing where $b$ takes large or small values provides information about where a future observation $x$ of $X$ will have greatest leverage on the value of $\int_{\mathcal{I}} bx$. Such information can be very useful for understanding the role played by the functional explanatory variable. Nevertheless, as this example suggests, the greatest overall interest lies, as in conventional linear regression, in using an estimator $\hat{b}$ as an aid to predicting, either qualitatively or quantitatively, a future value of $\int_{\mathcal{I}} bx$.

Thus, while there is extensive literature on properties of $\hat{b}$, for example on convergence rates of $\hat{b}$ to $b$ (see, e.g., [11, 13, 15, 20]), there is arguably a still greater need to understand the manner in which $\hat{b}$ should be constructed in order to optimize the prediction of $\int_{\mathcal{I}} bx$, or of $a + \int_{\mathcal{I}} bx$. This is the problem addressed in the present paper.

Estimation of $b$ is intrinsically an infinite-dimensional problem. Therefore, unlike slope estimation in conventional finite-dimensional regression, it involves smoothing or regularization. The smoothing step is used to reduce dimension, and the extent to which this should be done depends on the use to which the estimator of $b$ will be put, as well as on the smoothness of $b$. It is in this way that the problem of estimating $\int_{\mathcal{I}} bx$ is quite different from that of estimating $b$. The operation of integration, in computing $\int_{\mathcal{I}} \hat{b} x$ from $\hat{b}$, confers additional smoothness, with the result that if we smooth $\hat{b}$ optimally for estimating $b$ then it will usually be oversmoothed for estimating $\int_{\mathcal{I}} bx$.

Therefore the construction of $\hat{b}$, as a prelude to estimating $\int_{\mathcal{I}} bx$, should involve significant undersmoothing relative to the amount of smoothing that would be used if we wished only to estimate $b$ itself. In fact, as we shall show, the degree of undersmoothing can be so great that it enables $\int_{\mathcal{I}} bx$ to be estimated root-$n$ consistently, even though $b$ itself could not be estimated at such a fast rate.

However, root-$n$ consistency is not always possible when estimating $\int_{\mathcal{I}} bx$. The optimal convergence rate depends on a delicate balance among the smoothness of $b$, the smoothness of $x$, and the smoothness of the autocovariance of the stochastic process $X$, all measured with respect to the same sequence of basis functions. In a qualitative sense, $\int_{\mathcal{I}} bx$ can be estimated root-$n$ consistently if and only if $x$ is sufficiently smooth relative to the degree of smoothness of the autocovariance. If $x$ is less smooth than this, then the optimal rate at which $\int_{\mathcal{I}} bx$ can be estimated is determined jointly by the smoothnesses of $b$, $x$ and the autocovariance, and becomes faster as the smoothnesses of $x$ and of $b$ increase, and also as the smoothness of the covariance decreases.

These results are made explicitly clear in Section 4, which gives upper bounds to rates of convergence for specific estimators of $\int_{\mathcal{T}} bx$, and lower bounds (of the same order as the upper bounds) to rates of convergence for general estimators. Section 2 describes construction of the specific estimators of $b$, which are then substituted for $b$ in the formula $\int_{\mathcal{T}} bx$. Practical choice of smoothing parameters is discussed in Section 3.

In this brief account of the problem we have omitted mention of the role of the intercept, $a$, in the prediction problem. It turns out that from a theoretical viewpoint the role is minor. Given an estimator $\hat{b}$ of $b$, we can readily estimate $a$ by $\hat{a} = \bar{Y} - \int_{\mathcal{T}} \hat{b} \bar{X}$, where $\bar{X}$ and $\bar{Y}$ denote the means of the samples of $X_i$'s and $Y_i$'s, respectively. Taking this approach, it emerges that the rate of convergence of our estimator of $a + \int_{\mathcal{T}} bx$ is identical to that of our estimator of $\int_{\mathcal{T}} bx$, up to terms that converge to zero at the parametric rate $n^{-1/2}$. This point will be discussed in greater detail in Section 4.1.

The approach taken in this paper to estimating $b$ is based on functional principal components. While other methods could be used, the PC technique is currently the most popular. It goes back to work of Besse and Ramsay [1], Ramsay and Dalzell [27], Rice and Silverman [31] and Silverman [32, 33]. There are a great many more recent contributions, including those of Brumback and Rice [5], Cardot [7], Cardot, Ferraty and Sarda [8, 9, 10], Girard [19], James, Hastie and Sugar [23], Boente and Fraiman [3] and He, Müller and Wang [21].

Other recent work on regression for functional data includes that of Ferré and Yao [18], who introduced a functional version of sliced inverse regression; Preda and Saporta [26], who discussed linear regression on clusters of functional data; Escabias, Aguilera and Valderrama [14] and Ratcliffe, Heller and Leader [30], who described applications of functional logistic regression; and Ferraty and Vieu [16, 17] and Masry [24], who addressed various aspects of nonparametric regression for functional data. Müller and Stadtmüller [25] introduced the generalized functional linear model, where the response $Y_i$ is a general smooth function of $a + \int_{\mathcal{T}} bX_i$, plus an error. See also [22] and [12]. The methods developed in the present paper could be extended to this setting.

**2. Model and estimators.** We shall assume model (1.1), and suppose that the errors $\varepsilon_i$ are independent and identically distributed with zero mean and finite variance. It will be assumed too that the errors are independent of the $X_i$'s and that $\int_{\mathcal{T}} E(X^2) < \infty$.

Conventionally, estimation of $b$ is undertaken using a principal components approach, as follows. We take the covariance function of $X$ to be positive definite, in which case it admits a spectral decomposition in terms of strictly positive eigenvalues $\theta_j$,

$$(2.1) \quad K(u,v) \equiv \text{cov}\{X(u), X(v)\} = \sum_{j=1}^{\infty} \theta_j \phi_j(u) \phi_j(v), \qquad u, v \in \mathcal{I},$$

where $(\theta_j, \phi_j)$ are (eigenvalue, eigenfunction) pairs for the linear operator with kernel $K$, the eigenvalues are ordered so that $\theta_1 > \theta_2 > \cdots$ (in particular, we assume there are no ties among the eigenvalues), and the functions $\phi_1, \phi_2, \ldots$ form an orthonormal basis for the space of all square-integrable functions on $\mathcal{I}$.

Empirical versions of $K$ and of its spectral decomposition are

$$\widehat{K}(u,v) \equiv \frac{1}{n} \sum_{i=1}^{n} \{X_i(u) - \bar{X}(u)\}\{X_i(v) - \bar{X}(v)\}$$

$$= \sum_{j=1}^{\infty} \hat{\theta}_j \hat{\phi}_j(u) \hat{\phi}_j(v), \qquad u, v \in \mathcal{I},$$

where $\bar{X} = n^{-1} \sum_i X_i$. Analogously to the case of $K$, $(\hat{\theta}_j, \hat{\phi}_j)$ are (eigenvalue, eigenfunction) pairs for the linear operator with kernel $\widehat{K}$, ordered such that $\hat{\theta}_1 \geq \hat{\theta}_2 \geq \cdots$. Moreover, $\hat{\theta}_j = 0$ for $j \geq n + 1$. We take $(\hat{\theta}_j, \hat{\phi}_j)$ to be our estimator of $(\theta_j, \phi_j)$. The function $b$ can be expressed in terms of its Fourier series, as $b = \sum_{j \geq 1} b_j \phi_j$, where $b_j = \int b \phi_j$. We estimate $b$ as

$$(2.2) \qquad \hat{b} = \sum_{j=1}^{m} \hat{b}_j \hat{\phi}_j,$$

where $m$, lying in the range $1 \leq m \leq n$, denotes a "frequency cut-off" and $\hat{b}_j$ is an estimator of $b_j$.

To construct $\hat{b}_j$ we note that $b_j = \theta_j^{-1} g_j$, where $g_j$ denotes the $j$th Fourier coefficient of $g(u) = \int_{\mathcal{I}} K(u,v) b(v) \, dv$. A consistent estimator of $g$ is given by

$$\hat{g}(t) = \frac{1}{n} \sum_{i=1}^{n} \{X_i(t) - \bar{X}(t)\}(Y_i - \bar{Y}),$$

and so, for $1 \leq j \leq m$, we take $\hat{b}_j = \hat{\theta}_j^{-1} \hat{g}_j$, where $\hat{g}_j = \int_{\mathcal{I}} \hat{g} \hat{\phi}_j$.

While the problem of estimating $b$ is of intrinsic interest, it is arguably not as much practical importance as that of prediction, that is, estimating

$$p(x) \equiv E(Y|X = x) = a + \int_{\mathcal{I}} bx$$

for a particular function $x$. To accomplish this task we require an estimator of $a$,

$$\hat{a} = \bar{Y} - \int_{\mathcal{I}} \hat{b}\bar{X} = a - \int_{\mathcal{I}}(\hat{b} - b)\bar{X} + \bar{\varepsilon}.$$

Here, $\bar{Y}$ and $\bar{\varepsilon}$ are the respective means of the sequences $Y_i$ and $\varepsilon_i$. Our estimator of $p(x)$, for a given function $x$, is

$$\hat{p}(x) = \hat{a} + \int_{\mathcal{I}} \hat{b}x.$$

In Section 4 we shall introduce three parameters, $\alpha$, $\beta$ and $\gamma$, describing the smoothness of $K$, $b$ and $x$, respectively. In each case, smoothness is measured in the context of generalized Fourier expansions in the basis $\phi_1, \phi_2, \ldots$, and the larger the value of the parameter, the smoother the associated function. We shall show in Theorem 4.1 that if $x$ is sufficiently smooth relative to $K$, specifically if $\gamma > \frac{1}{2}(\alpha + 1)$, then $\int_{\mathcal{I}} bx$ can be estimated root-$n$ consistently. For smaller values of $\gamma$, the optimal convergence rate is slower than $n^{-1/2}$.

**3. Numerical implementation and simulation study.** There is a variety of possible approaches to empirical choice of the cut-off, $m$, although not all are directly suited to estimation of $\int_{\mathcal{I}} bx$. Potential methods include those based on simple least-squares, on the bootstrap or on cross-validation. In some instances where $\int_{\mathcal{I}} \hat{b}x$ is root-$n$ consistent for $\int_{\mathcal{I}} bx$, $m$ can be chosen within a wide range without appreciably affecting the performance of the estimator. Only in relatively "unsmooth" cases, where either $\gamma \leq \frac{1}{2}(\alpha + 1)$, or $\gamma > \frac{1}{2}(\alpha + 1)$ but $\gamma$ is close to $\frac{1}{2}(\alpha + 1)$, is the choice of $m$ rather critical. The empirical identification of unsmooth cases, and empirical choice of $m$ in those instances, are challenging problems, and we shall not attempt to address them here. (See the last paragraph of Section 2 for discussion of $\alpha$, $\beta$ and $\gamma$.)

Instead, we shall give below a simple threshold-based algorithm for choosing $m$ empirically in cases where $x$ is sufficiently smooth. There, the algorithm guarantees root-$n$ consistency. The order of magnitude of the empirically chosen $m$ depends very much on selection of the threshold, but nevertheless the estimator $\int_{\mathcal{I}} \hat{b}x$ remains root-$n$ consistent in a very wide range of cases. Therefore, the effectiveness of the threshold algorithm underscores the robustness of the estimator against choice of $m$ in cases where $x$ is smooth.

To describe the threshold algorithm, let $C > 0$ and $0 < c \leq \frac{1}{2}$, and put $I_j = 1$ if $\hat{\theta}_j \geq t \equiv Cn^{-c}$, with $I_j = 0$ otherwise. Since the sequence $\hat{\theta}_1, \hat{\theta}_2, \ldots$ is nonincreasing and $\hat{\theta}_j = 0$ for $j \geq n + 1$, then $I_1, I_2, \ldots$ is a sequence of $\hat{m}$, say, 1's, followed by an infinite sequence of 0's. Therefore the threshold algorithm implicitly gives an empirical rule for choosing the cut-off, $m$. Our estimator of $\int_{\mathcal{I}} bx$ is $\int_{\mathcal{I}} \hat{b}x$, where $\hat{b} = \sum_{1 \leq j \leq \hat{m}} \hat{b}_j \hat{\phi}_j$. Note that the estimator

$$\int_{\mathcal{I}} \hat{b}x = \sum_j I_j \hat{b}_j \bar{x}_j = \sum_{1 \leq j \leq \hat{m}} \hat{b}_j \bar{x}_j,$$

where $\bar{x}_j = \int_{\mathcal{I}} x \hat{\phi}_j$. This form is often easier to use in numerical calculations.

To appreciate the size of $\hat{m}$ chosen by this rule, let us suppose that $\theta_j = \text{const}.j^{-\alpha}$. It can be shown that, for the specified range of values of $c$, $\hat{\theta}_j = \text{const}.j^{-\alpha}\{1 + o_p(1)\}$ uniformly in $1 \leq j \leq \hat{m} + k$, for each integer $k \geq 1$. Therefore, $\hat{m} = \text{const}.n^{c/\alpha}\{1 + o_p(1)\}$. It follows that the order of magnitude of $\hat{m}$ changes a great deal as we vary $c$.

It can be proved too that, under the conditions of Theorem 4.1, and assuming that $\alpha \geq 2$, $\gamma \geq \frac{3}{2}(\alpha + 2)$ and $\beta + \gamma \geq (\alpha/2c) + 1$,

$$(3.1) \qquad \sum_{j=1}^{\hat{m}} \hat{b}_j \bar{x}_j = \int_{\mathcal{I}} bx + O_p(n^{-1/2}).$$

This result demonstrates the root-$n$ consistency of the estimator on the left-hand side, for a range of different orders of magnitude of $\hat{m}$. Of course, (3.1) continues to hold if the number of terms, $\hat{m}$, is replaced by a deterministic quantity, say $m \sim \text{const}.n^{c/\alpha}$. Note too that the conditions $\gamma \geq \frac{3}{2}(\alpha + 2)$ and $\beta + \gamma \geq (\alpha/2c) + 1$ are both implied by $\gamma \geq \max(3/2, 1/2c)\alpha + 3$, which asserts simply that the function $x$ is sufficiently smooth relative to $K$.

The case where the functions $X_i$ are observed on a regular grid of $k$ points with additive white noise may be treated similarly. Indeed, it can be proved that if continuous approximations to the $X_i$'s are generated by passing a local-linear smoother through noisy, gridded data, and if we take $c = \frac{1}{2}$, then all the results discussed above remain true provided $n = O(k)$. That is, $k$ should be of the same order as, or of larger order than, $n$. Details are given in the Appendix of [6]. Similar results are obtained using smoothing methods based on splines or orthogonal series.

A simulation study was carried out to investigate the finite-sample performance of the thresholding procedure given above. The study considered the model (1.1) in two cases. In the first, the predictor $X_i$ was observed continuously without error. Specifically, random samples of size $n = 100$ were generated from the model (1.1), where $\mathcal{I} = [0, 1]$, the random functions $X_i$ were distributed as $X = \sum_j Z_j 2^{1/2} \cos(j\pi t)$, the $Z_j$'s were independent and normal $N(0, 4j^{-2})$, $b = \sum_j j^{-4} 2^{1/2} \cos(j\pi t)$, and the errors $\varepsilon_i$ were independent and normal $N(0, 4)$. The future observation of $X$ was taken to be $x = \sum_j j^{-2} 2^{1/2} \cos(j\pi t)$, in which case the conditional mean of $y$ given $X = x$ was 1.0141.

TABLE 1
*Comparison of average squared errors*

| Threshold | 0.001 | 0.01 | 0.05 | 0.1 | 0.15 | 0.2 |
|---|---|---|---|---|---|---|
| X continuous | 0.026 | 0.019 | 0.015 | 0.014 | 0.013 | 0.015 |
| X discrete with noise | 0.035 | 0.022 | 0.016 | 0.017 | 0.015 | 0.016 |

The example in the second case was the same as that for the first, except that each $X_i$ was observed discretely on an equally-spaced grid of 200 points with additive $N(0, 1)$ random noise. We used an orthogonal-series smoother to "estimate" each $X_i$ from the corresponding discrete data. Table 1 gives values of averaged squared error of the estimator of the conditional mean, computed by averaging 500 Monte Carlo simulations. It is clear from these results that the procedure is robust against discretization, random errors and choice of the threshold.

Earlier in this section we discussed the robustness of $\hat{b}$ to choice of smoothing parameter in the prediction problem. This robustness is not shared in cases where $\hat{b}$ is of interest in its own right, rather than a tool for prediction. To make this comparison explicit, and to compare the levels of smoothing appropriate for prediction and estimation, we extended the simulation study above. We selected $X$ as before, but took $b = 10 \sum_j j^{-2} 2^{1/2} \cos(j\pi t)$ and $x = \sum_j j^{-1.6} 2^{1/2} \cos(j\pi t)$. In the case of noisy, discrete observations we took the noise to be $N(0, 1)$ and the grid to consist of 500 points. Sample size was $n = 100$.

For the thresholds $t = 0.001, 0.01, 0.05, 0.1, 0.15, 0.2$ used to construct Table 1, mean squared prediction error was relatively constant; respective values were $0.013, 0.008, 0.007, 0.010, 0.015, 0.022$. However, mean integrated squared error of $\hat{b}$ was as high as 168 when $t = 0.001$, dropping to 6.67 at $t = 0.01$ and reaching its minimum, 0.639, at $t = 0.1$. Similar results were achieved in the case of noisy, discrete data; values of mean squared prediction error there were $0.014, 0.008, 0.009, 0.013, 0.019, 0.028$ for the respective values of $t$, and mean integrated squared error of $\hat{b}$ was elevated by about 30% across the range, the minimum again occurring when $t = 0.1$.

These results also indicate the advantages of undersmoothing when making predictions, as opposed to estimating $\hat{b}$ in its own right. In particular, the numerical value of the optimal threshold for prediction is a little less than that for estimating $\hat{b}$. Discussion of theoretical aspects of this point will be given in Section 4.

## 4. Convergence rates.

4.1. *Effect of the intercept, a.* In terms of convergence rates, the problems of estimating $a + \int_{\mathcal{I}} bx$ and $\int_{\mathcal{I}} bx$ are not intrinsically different. To appreciate this point, define $\mu = E(X)$, let the functionals $p$ and $\hat{p}$ be as in Section 2, and put $q(x) = \int_{\mathcal{I}} b(x - \mu)$ and $\hat{q}(x) = \int_{\mathcal{I}} \hat{b}(x - \mu)$. Given a random variable $Z$, write $M(Z) = (EZ^2)^{1/2}$. Then

$$
\begin{aligned}
|M\{\hat{p}(x) &- p(x)\} - M\{\hat{q}(x) - q(x)\}| \\
&\leq M\left\{\int_{\mathcal{I}} (\hat{b} - b)(\bar{X} - \mu) + \bar{\varepsilon}\right\} \\
&\leq (E\|\hat{b} - b\|^2)^{1/2}(E\|\bar{X} - \mu\|^2)^{1/2} + (E\bar{\varepsilon}^2)^{1/2}.
\end{aligned}
\tag{4.1}
$$

Provided only that $E\|\hat{b} - b\|^2$ is bounded, the right-hand side of (4.1) equals $O(n^{-1/2})$. Hence, (4.1) shows that, up to terms that converge to zero at the parametric rate $n^{-1/2}$, the rates of convergence of $\hat{p}(x)$ to $p(x)$ and of $\hat{q}(x)$ to $q(x)$ are identical. This result, and the fact that $q(x)$ is identical to $\int bx$ provided $x$ is replaced by $x - \mu$, imply that when addressing convergence rates in the prediction problem it is sufficient to treat estimation of $\int_{\mathcal{I}} bx$.

4.2. *Estimation of $\int bx$.*  Recall that our estimator of $\int bx$ is $\int \hat{b}x$. Suppose the eigenvalues $\theta_j$ in the spectral decomposition (2.1) satisfy

$$
C^{-1}j^{-\alpha} \leq \theta_j \leq Cj^{-\alpha}, \qquad \theta_j - \theta_{j+1} \geq C^{-1}j^{-\alpha-1} \qquad \text{for } j \geq 1.
\tag{4.2}
$$

For example, if $\theta_j = Dj^{-\alpha}$ for a constant $D > 0$, then $\theta_j - \theta_{j+1} \sim D\alpha^{-1}j^{-\alpha-1}$, and so (4.2) holds. The second part of (4.2) asks that the spacings among eigenvalues not be too small. Methods based on a frequency cut-off $m$ can have difficulty when spacings equal zero, or are close to zero. To appreciate why, note that if $\theta_{j+1} = \cdots = \theta_{j+k}$ then $\phi_{j+1}, \ldots, \phi_{j+k}$ are not individually identifiable (although the set of these $k$ functions is identifiable). In particular, individual functions cannot be estimated consistently. This can cause problems when estimating $\int_{\mathcal{I}} bx$ if the frequency cut-off lies strictly between $j$ and $j + k$.

Let $Z$ have the distribution of a generic $X_i - E(X_i)$. Then we may write $Z = \sum_{j \geq 1} \xi_j \phi_j$, where $\xi_j = \int Z \phi_j$ is the $j$th principal component, or Karhunen–Loève coefficient, of $Z$. We assume that all the moments of $X$ are finite, and more specifically that for each $r \geq 2$ and each $j \geq 1$, $E|\xi_j|^{2r} \leq C(r)\theta_j^r$, where $C(r)$ does not depend on $j$; and, for any sequence $j_1, \ldots, j_4$, $E(\xi_{j_1} \ldots \xi_{j_4}) = 0$ unless each index $j_k$ is repeated.
$\tag{4.3}$

In particular, (4.3) holds if $X$ is a Gaussian process. Let $\beta > 1$ and $C_1 > 0$, and let

$$
\mathcal{B} = \mathcal{B}(C_1, \beta) = \left\{b : b = \sum_{j \geq 1} b_j \phi_j, \text{with } |b_j| \leq C_1 j^{-\beta} \text{ for each } j \geq 1\right\}.
\tag{4.4}
$$

We can interpret $\mathcal{B}(C_1,\beta)$ as a "smoothness class" of functions, where the functions become smoother (measured in the sense of generalized Fourier expansions in the basis $\phi_1,\phi_2,\dots$) as $\beta$ increases. We suppose too that the fixed function $x$ satisfies

$$(4.5) \qquad x = \sum_{j=1}^{\infty} x_j \phi_j \qquad \text{with } |x_j| \le C_2 j^{-\gamma} \text{ for each } j.$$

Again, $x$ becomes smoother in the sense of generalized Fourier expansions as $\gamma$ increases.

Define $m_0 = m_0(n)$ by

$$(4.6) \qquad m_0 = \begin{cases} n^{1/2(\beta+\gamma-1)}, & \text{if } \alpha+1 < 2\gamma, \\ (n/\log n)^{1/(\alpha+2\beta-1)}, & \text{if } \alpha+1 = 2\gamma, \\ n^{1/(\alpha+2\beta-1)}, & \text{if } \alpha+1 > 2\gamma. \end{cases}$$

These explicit values serve to simplify our discussion and our proof of Theorem 4.1, and do not reflect the wider range of values of $m$, particularly in the case $\alpha+1 < 2\gamma$, for which our theory is valid. Discussion of this point has been given in Section 3.

Recall the definition of $\hat{b}$ at (2.2). Given arbitrary positive constants $C_3$, $C_4$ and $C_5$, let

$$(4.7) \qquad \tilde{b} = \begin{cases} \hat{b}, & \text{if } \|\hat{b}\| \le C_4 n^{C_5}, \\ C_3, & \text{otherwise}, \end{cases}$$

where, for a function $\psi$ on $\mathcal{I}$, $\|\psi\|^2 = \int_{\mathcal{I}} \psi^2$. This truncation of $\hat{b}$ serves to ensure that all moments of $\tilde{b}$ are finite.

THEOREM 4.1. *Assume the eigenvalues $\theta_j$ satisfy (4.2), that (4.3) holds and that all moments of the distribution of the errors $\varepsilon_i$ are finite. Let $\alpha$, $\beta$ and $\gamma$ be as in (4.2), (4.4) and (4.5), respectively. Suppose that $\alpha > 1$, $\beta \ge \alpha + 2$ and $\gamma > \frac{1}{2}$, and that the ratio of $m$ to $m_0$ is bounded away from zero and infinity as $n \to \infty$. Then, for each given $C, C_1, \dots, C_5 > 0$, as $n \to \infty$, the estimator $\tilde{b}$ given in (4.7) satisfies*

$$(4.8) \qquad \sup_{b \in \mathcal{B}(C_1,\beta)} E\left( \int_{\mathcal{I}} \tilde{b} x - \int_{\mathcal{I}} b x \right)^2 = O(\tau),$$

*where $\tau = \tau(n)$ is given by*

$$(4.9) \qquad \tau = \begin{cases} n^{-1}, & \text{if } \alpha+1 < 2\gamma, \\ n^{-1}\log n, & \text{if } \alpha+1 = 2\gamma, \\ n^{-2(\beta+\gamma-1)/(\alpha+2\beta-1)}, & \text{if } \alpha+1 > 2\gamma. \end{cases}$$

The smoothing-parameter choices suggested by (4.6) are different from those that would be used if our aim were to estimate $b$ rather than $\int_{\mathcal{I}} bx$. In particular, to optimize the $L_2$ convergence rate of $\tilde{b}$ to $b$ we would take $m$ to be of size $n^{1/(\alpha+2\beta)}$ in each of the three settings addressed by (4.6). See, for example, [20]. In the critical cases where $\alpha + 1 \geq 2\gamma$, this provides an order of magnitude more smoothing than is suggested by (4.6). The intuition behind this result is that the integration step, in the definition $\int_{\mathcal{I}} \hat{b}x$, provides additional smoothing no matter what level is used when constructing $\hat{b}$, and so less smoothing is needed for $\hat{b}$.

The case $\alpha + 1 < 2\gamma$ is more difficult to discuss in these terms, since a variety of different orders of magnitude of $m$ can lead to the same optimal mean-square convergence rate of $n^{-1}$. Further discussion of this issue is given in Section 3.

Of course, there are other related problems where similar phenomena are observed. Consider, for example, the problem of estimating a distribution function by integrating a kernel density estimator. In order to achieve the same parametric convergence rate as the empirical distribution function, we should, when constructing the density estimator, use a substantially smaller bandwidth than would be appropriate if we wanted a good estimator of the density itself. The operation of integrating the density estimator provides additional smoothing, over and above that accorded by the bandwidth, and so if the net result is not to be an oversmoothed distribution-function estimator then we should smooth less at the density estimation step. The same is true in the problem of prediction in functional regression; the operation of integrating $\tilde{b}x$ provides additional smoothing, and so to get the right amount of smoothing in the end we should undersmooth when computing the slope-function estimator. A curious feature of the regression prediction problem is that, unlike the distribution estimation one, it is not always parametric, and in some cases the optimal convergence rate lies strictly between that for the nonparametric problem of slope estimation and the parametric $n^{-1/2}$ rate.

4.3. *Lower bounds.* We adopt notation from Sections 4.1 and 4.2, and in particular take $x = \sum_{j \geq 1} x_j \phi_j$ to be a function and define $\mathcal{B}$ as at (4.4). Recall that the functions $\phi_j$ form an orthonormal basis for square-integrable functions on $\mathcal{I}$. Assume that, for a constant $C_6 > 1$,

$$C_6^{-1} \leq j^\alpha \theta_j \leq C_6 \quad \text{and} \quad C_6^{-1} \leq j^\gamma |x_j| \leq C_6 \qquad \text{for all } j \geq 1.$$

Let $\widehat{T}$ denote any estimator of $T(b) = \int_{\mathcal{I}} bx$, and define $\tau = \tau(n)$ as at (4.9). Our main result in this section provides a lower bound to the convergence rate of $\widehat{T}$ to $T(b)$, complementing the upper bound given by Theorem 4.1 in the case $\widehat{T} = \int_{\mathcal{I}} \tilde{b}x$, where $\tilde{b}$ is given by (4.7). We make relatively specific assumptions about the nature of the model, for example that $X$ is a Gaussian process and the intercept, $a$, vanishes, bearing in mind that in the case of a lower bound, the strength of the result is increased, from some viewpoints, through imposing relatively narrow conditions.

THEOREM 4.2. *Let $\alpha$, $\beta$ and $\gamma$ be as in (4.2), (4.4) and (4.5), respectively, and assume $\alpha, \beta > 1$ and $\gamma > \frac{1}{2}$. Suppose too that the process $X$ is Gaussian and that the errors $\varepsilon_i$ in the model (1.1) are Normal with zero mean and strictly positive variance; and take $a = 0$. Then there exists a constant $C_7 > 0$ such that, for any estimator $\hat{T}$ and for all sufficiently large $n$,*

$$\sup_{b \in \mathcal{B}(C_1, \beta)} E\{\hat{T} - T(b)\}^2 \geq C_7 \tau,$$

*where $\tau = \tau(n)$ is given as in (4.9).*

A comparison of the lower bound given above with the upper bound given in Theorem 4.1 yields the result that the minimax risk of estimating $\int bx$ satisfies

$$\inf_{\hat{T}} \sup_{b \in \mathcal{B}(C_1, \beta)} E\left(\hat{T} - \int bx\right)^2 \asymp \begin{cases} n^{-1}, & \text{if } \alpha + 1 < 2\gamma, \\ n^{-1} \log n, & \text{if } \alpha + 1 = 2\gamma, \\ n^{-2(\beta+\gamma-1)/(\alpha+2\beta-1)}, & \text{if } \alpha + 1 > 2\gamma, \end{cases}$$

where, for positive sequences $a_n$ and $b_n$, $a_n \asymp b_n$ means that $a_n / b_n$ is bounded away from zero and infinity as $n \to \infty$.

## 5. Proof of Theorem 4.1.

5.1. *Preliminaries.* Define $\Delta = \widehat{K} - K$, $\|\Delta\|^2 = \int_{\mathcal{I}^2} \Delta^2$ and $\delta_j = \min_{k \leq j}(\theta_k - \theta_{k+1})$. It may be shown from results of Bhatia, Davis and McIntosh [2] that

$$(5.1) \qquad \begin{aligned} \sup_{j \geq 1} |\hat{\theta}_j - \theta_j| &\leq \|\Delta\|, \\ \sup_{j \geq 1} \delta_j \|\hat{\phi}_j - \phi_j\| &\leq 8^{1/2} \|\Delta\|. \end{aligned}$$

For simplicity in our proof we shall take $m = m_0$, as defined in (4.6). Note that in this setting $m \leq n^{1/(\alpha+2\beta-1)}$ in each of the three cases in (4.6).

Expand $x$ with respect to both the orthonormal series $\phi_1, \phi_2, \ldots$ and $\hat{\phi}_1, \hat{\phi}_2, \ldots$, obtaining $x = \sum_{j \geq 1} x_j \phi_j = \sum_{j \geq 1} \bar{x}_j \hat{\phi}_j$, where $x_j = \int_{\mathcal{I}} x \phi_j$ and $\bar{x}_j = \int_{\mathcal{I}} x \hat{\phi}_j$. Put $\tilde{g}_j = \int_{\mathcal{I}} g \phi_j$. In this notation

$$\int_{\mathcal{I}} (\hat{b} - b)x = \sum_{j=1}^{m} (\hat{b}_j \bar{x}_j - b_j x_j) - \sum_{j=m+1}^{\infty} b_j x_j,$$

whence it follows that

$$
\begin{aligned}
\left| \int_{\mathcal{I}} (\hat{b} - b) x \right| &\leq \left| \sum_{j=1}^{m} (\hat{b}_j - b_j) x_j \right| \\
&\quad + \left| \sum_{j=m+1}^{\infty} b_j x_j \right| + \left| \sum_{j=1}^{m} b_j (\bar{x}_j - x_j) \right| \\
&\quad + \sum_{j=1}^{m} |\hat{b}_j - b_j| |\bar{x}_j - x_j|.
\end{aligned}
\tag{5.2}
$$

It is straightforward to show that $|\sum_{j \geq m+1} b_j x_j| = O(m^{-(\beta+\gamma-1)})$. This quantity equals $O\{(n^{-1} \log n)^{1/2}\}$ if $\alpha + 1 = 2\gamma$, equals $O(n^{-(\beta+\gamma-1)/(\alpha+2\beta-1)})$ if $\alpha + 1 > 2\gamma$ and equals $o(n^{-1/2})$ otherwise. We shall complete the derivation of Theorem 4.1 by obtaining bounds for second moments of the other three terms on the right-hand side of (5.2). Our analysis will show that the first and second terms determine the convergence rate, and that the third and fourth terms are asymptotically negligible. In the arguments leading to the bounds we shall use the notation "const." to denote a constant, the value of which does not depend on $b \in \mathcal{B}$. In particular, the bounds we shall give are valid uniformly in $b$, although we shall not mention that property explicitly.

5.2. *Bound for* $|\sum_{j \leq m} (\hat{b}_j - b_j) x_j|$. Note that

$$
\hat{b}_j - b_j = (\hat{\theta}_j^{-1} - \theta_j^{-1})(\hat{g}_j - g_j) + \theta_j^{-1}(\hat{g}_j - g_j) + (\hat{\theta}_j^{-1} - \theta_j^{-1}) g_j,
\tag{5.3}
$$

$$
\hat{g}_j - g_j = \tilde{g}_j - g_j + \int_{\mathcal{I}} (\hat{g} - g)(\hat{\phi}_j - \phi_j) + \int_{\mathcal{I}} g(\hat{\phi}_j - \phi_j).
\tag{5.4}
$$

Therefore, defining $\Delta_g = \hat{g} - g$, we have

$$
\left| \hat{g}_j - g_j - \int_{\mathcal{I}} g(\hat{\phi}_j - \phi_j) \right| \leq 3 \|\Delta_g\|.
\tag{5.5}
$$

If the event

$$
\mathcal{E} = \{ |\hat{\theta}_j - \theta_j| \leq \tfrac{1}{2} \theta_j \text{ for all } 1 \leq j \leq m \}
\tag{5.6}
$$

holds, then $|\hat{\theta}_j^{-1} - \theta_j^{-1}| \leq 2|\hat{\theta}_j - \theta_j|/\theta_j^2 \leq \theta_j^{-1}$. It can be proved, using this result, (5.1), (5.4) and (5.5), that if $\mathcal{E}$ holds,

$$
\begin{aligned}
\tfrac{1}{6} \left| \sum_{j=1}^{m} (\hat{b}_j - b_j) x_j \right| &\leq \left| \sum_{j=1}^{m} (\tilde{g}_j - g_j) x_j \theta_j^{-1} \right| + \left| \sum_{j=1}^{m} x_j \theta_j^{-1} \int_{\mathcal{I}} g(\hat{\phi}_j - \phi_j) \right| \\
&\quad + \|\Delta\| \sum_{j=1}^{m} \left| \int_{\mathcal{I}} g(\hat{\phi}_j - \phi_j) \right| |x_j| \theta_j^{-2}
\end{aligned}
\tag{5.7}
$$

$$+ 8^{1/2}\|\Delta\|\sum_{j=1}^{m}(\|\Delta_g\|\delta_j^{-1} + |g_j|\theta_j^{-1})|x_j|\theta_j^{-1}.$$

For each real number $r$, define

$$t_r(m) = \begin{cases} m^{r+1}, & \text{if } r > -1, \\ \log m, & \text{if } r = -1, \\ 1, & \text{if } r < -1. \end{cases}$$

Standard moment calculations, noting that $S_1(g) \equiv \sum_{j \leq m}(\tilde{g}_j - g_j)x_j\theta_j^{-1}$ may be expressed as a sum of $n$ independent and identically distributed random variables with zero mean, show that $E\{S_1(g)^2\} \leq \text{const.}\,n^{-1}t_{2\alpha-2\gamma}(m)$, uniformly in $g$. Moreover, denoting by $S_2(g)$ the last term on the right-hand side of (5.7), we deduce that

$$
\begin{aligned}
E\{S_2(g)^2\} &\equiv E\Bigg\{\|\Delta\|\sum_{j=1}^{m}(\|\Delta_g\|\delta_j^{-1} + |g_j|\theta_j^{-1})|x_j|\theta_j^{-1}\Bigg\}^2 \\
&\leq \text{const.}\{n^{-2}t_{2\alpha-\gamma+1}(m)^2 + n^{-1}t_{\alpha-\beta-\gamma}(m)^2\}.
\end{aligned}
$$
(5.8)

If $\beta \geq \gamma$ then $t_{\alpha-\beta-\gamma}(m) \leq t_{\alpha-2\gamma}(m)$, and if $\beta < \gamma$ then, since $\beta > \frac{1}{2}(\alpha+1)$, $\alpha - \beta - \gamma < -1$, implying that $t_{\alpha-\beta-\gamma}(m) \leq \text{const.}\,t_{\alpha-2\gamma}(m)$. Moreover, $t_{2\alpha-\gamma+1}(m) \leq \text{const.}\,t_{\alpha-2\gamma}(m)m^{\alpha+1}$, and by assumption, $n \geq m^{\alpha+1}$. Therefore, $n^{-1}t_{2\alpha-\gamma+1}(m) \leq \text{const.}\,t_{\alpha-2\gamma}(m)$. Hence, (5.8) implies that $E\{S_2(g)^2\} \leq \text{const.}\,n^{-1}t_{\alpha-2\gamma}(m)$. Combining this bound with that for $E\{S_1(g)^2\}$, and with (5.7), and writing $I(\mathcal{F})$ for the indicator function of any subset $\mathcal{F} \subseteq \mathcal{E}$, we deduce that

$$
\begin{aligned}
&E\Bigg[I(\mathcal{F})\Bigg\{\sum_{j=1}^{m}(\hat{b}_j - b_j)x_j\Bigg\}^2\Bigg] \\
&\leq \text{const.}\Bigg(E\Bigg[I(\mathcal{F})\Bigg\{\sum_{j=1}^{m}x_j\theta_j^{-1}\int_{\mathcal{I}}g(\hat{\phi}_j - \phi_j)\Bigg\}^2\Bigg] \\
&\qquad + E\Bigg[I(\mathcal{F})\|\Delta\|^2\Bigg\{\sum_{j=1}^{m}\Big|\int_{\mathcal{I}}g(\hat{\phi}_j - \phi_j)\Big||x_j|\theta_j^{-2}\Bigg\}^2\Bigg] \\
&\qquad\qquad\qquad\qquad\qquad + n^{-1}t_{\alpha-2\gamma}(m)\Bigg).
\end{aligned}
$$
(5.9)

Note too that if $\mathcal{E}$ holds,

$$
\begin{aligned}
\sum_{j=1}^{m}(\hat{b}_j - b_j)^2 &\leq \text{const.}\sum_{j=1}^{m}\theta_j^{-2}\Bigg\{(\tilde{g}_j - g_j)^2 + \Big|\int_{\mathcal{I}}g(\hat{\phi}_j - \phi_j)\Big|^2\Bigg\} \\
&\qquad + \text{const.}\|\Delta\|^2\{\|\Delta_g\|^2 t_{4\alpha+2}(m) + t_{2\alpha-2\beta}(m)\},
\end{aligned}
$$
(5.10)

and also that

$$
\left| \sum_{j=1}^{m} b_j (\bar{x}_j - x_j) \right| = \left| \sum_{j=1}^{m} b_j \int x (\hat{\phi}_j - \phi_j) \right|, \tag{5.11}
$$

$$
\sum_{j=1}^{m} (\bar{x}_j - x_j)^2 = \sum_{j=1}^{m} \left\{ \int_{\mathcal{I}} x (\hat{\phi}_j - \phi_j) \right\}^2. \tag{5.12}
$$

Let $p = g$ or $x$, and define $\pi = \alpha + \beta$ and $\pi = \gamma$ in the respective cases. Let $q_1, q_2, \ldots$ denote constants satisfying $|q_j| \le \text{const.} j^\kappa$ for each $j$, where $\kappa = \alpha - \gamma$ if $p = g$, and $\kappa = -(\alpha + \beta)$ if $p = x$. Given $\eta > 0$, consider the event

$$
\mathcal{F} = \{ \|\!|\Delta|\!\| \le n^{\eta - (1/2)} \text{ and } \tag{5.13}
$$
$$
|\hat{\theta}_j - \theta_j| \le \tfrac{1}{2} C^{-1} j^{-\alpha - 1} \text{ for all } 1 \le j \le m \}.
$$

Comparing (5.6) and (5.13), and noting (4.2), we see that $\mathcal{F} \subseteq \mathcal{E}$. We shall show in Section 5.5 that, uniformly in $1 \le j \le \text{const.} n^{1/(\alpha+1)}$,

$$
E \left\{ I(\mathcal{F}) \int_{\mathcal{I}} p (\hat{\phi}_j - \phi_j) \right\}^2 \le \text{const.} n^{-1} j^{-\alpha} (1 + j^{2\alpha + 2 - 2\pi}), \tag{5.14}
$$

and also,

$$
E \left[ I(\mathcal{F}) \left\{ \sum_{j=1}^{m} q_j \int_{\mathcal{I}} p (\hat{\phi}_j - \phi_j) \right\}^2 \right] \le \text{const.} n^{-1} t_{2\kappa - \alpha}(m). \tag{5.15}
$$

Next we use (5.15) to bound the first term on the right-hand side of (5.9):

$$
E \left[ I(\mathcal{F}) \left\{ \sum_{j=1}^{m} x_j \theta_j^{-1} \int_{\mathcal{I}} g (\hat{\phi}_j - \phi_j) \right\}^2 \right] \le \text{const.} n^{-1} t_{\alpha - 2\gamma}(m). \tag{5.16}
$$

To bound the second term, it can be proved from (5.14) that

$$
E \left[ I(\mathcal{F}) \left\{ \sum_{j=1}^{m} \left| \int_{\mathcal{I}} g (\hat{\phi}_j - \phi_j) \right| |x_j| \theta_j^{-2} \right\}^2 \right]
$$
$$
\le \text{const.} n^{-2\{\beta - \alpha - (3/2)\}/(\alpha + 2\beta - 1)}. \tag{5.17}
$$

Going back to the definition of $\mathcal{F}$ at (5.13), and taking $\eta < \{\beta - \alpha - (3/2)\}/(\alpha + 2\beta - 1)$, we deduce from (5.17) that

$$
E \left[ I(\mathcal{F}) \|\!|\Delta|\!\|^2 \left\{ \sum_{j=1}^{m} \left| \int_{\mathcal{I}} g (\hat{\phi}_j - \phi_j) \right| |x_j| \theta_j^{-2} \right\}^2 \right] \le \text{const.} n^{-1}. \tag{5.18}
$$

Results (5.9), (5.16) and (5.18) imply that

$$
E \left[ I(\mathcal{F}) \left\{ \sum_{j=1}^{m} (\hat{b}_j - b_j) x_j \right\}^2 \right] \le \text{const.} n^{-1} t_{\alpha - 2\gamma}(m). \tag{5.19}
$$

5.3. *Bounds for* $|\sum_{j\leq m} b_j(\bar{x}_j - x_j)|$ *and* $\sum_{j\leq m} |\hat{b}_j - b_j||\bar{x}_j - x_j|$. Noting that $\kappa = -(\alpha + \beta)$ when $p = x$, we may also use (5.15) and (5.14) to bound the expected values of the squares of the right-hand sides of (5.11) and (5.12), respectively, multiplied by $I(\mathcal{F})$:

$$(5.20) \qquad E\left[I(\mathcal{F})\left\{\sum_{j=1}^{m} b_j \int x(\hat{\phi}_j - \phi_j)\right\}^2\right] \leq \text{const.} n^{-1},$$

$$(5.21) \qquad E\left[I(\mathcal{F})\sum_{j=1}^{m}\left\{\int_{\mathcal{I}} x(\hat{\phi}_j - \phi_j)\right\}^2\right] \leq \text{const.} n^{-1} t_{\alpha+3-2\gamma}(m).$$

Noting that $\beta \geq \alpha + 2$ and $E(\tilde{g}_j - g_j)^2 \leq \text{const.} n^{-1}\theta_j$, we can show from (5.10) and (5.14) that

$$(5.22) \qquad E\left\{I(\mathcal{F})\sum_{j=1}^{m}(\hat{b}_j - b_j)^2\right\} \leq \text{const.} n^{-1} m^{\alpha+1}.$$

From (5.21) and (5.22) it follows that

$$
\begin{aligned}
&E\left[I(\mathcal{F})\left(\sum_{j=1}^{m}|\hat{b}_j - b_j||\bar{x}_j - x_j|\right)^2\right] \\
(5.23) \qquad &\leq E\left\{I(\mathcal{F})\sum_{j=1}^{m}(\hat{b}_j - b_j)^2\right\} E\left\{I(\mathcal{F})\sum_{j=1}^{m}(\bar{x}_j - x_j)^2\right\} \\
&\leq \text{const.} n^{-1} m^{\alpha+1} \cdot n^{-1} t_{\alpha+3-2\gamma}(m) \leq \text{const.} n^{-1}.
\end{aligned}
$$

5.4. *Completion of the proof of Theorem* 4.1. Combining (5.2), (5.19), (5.20) and (5.23) we deduce that

$$(5.24) \qquad E\left[I(\mathcal{F})\left\{\int_{\mathcal{I}}(\hat{b} - b)x\right\}^2\right] \leq \text{const.} n^{-1} t_{\alpha-2\gamma}(m).$$

The proof of Theorem 4.1 will be complete if we show that the factor $I(\mathcal{F})$ can be removed from the left-hand side. Since, in view of (4.7), our estimator $\tilde{b}$ satisfies $\|\tilde{b}\| \leq C_4 n^{C_5}$, then it suffices to prove that, for all $D > 0$, $P(\mathcal{F}) = 1 - O(n^{-D})$. Now the first part of (5.1) and (5.13) imply that if we define

$$\mathcal{G} = \{\|\Delta\| \leq \min(n^{\eta - (1/2)}, cC^{-1}m^{-\alpha-1})\},$$

then $\mathcal{G} \subseteq \mathcal{F}$. Since $m \leq n^{1/(\alpha+2\beta-1)}$ and $2(\alpha+1) < \alpha + 2\beta - 1$, then for some $\eta' > 0$, $m^{-\alpha-1} \geq n^{\eta'-(1/2)}$. Therefore, if $\zeta > 0$ is sufficiently small, there exists $n_0 \geq 1$ such that, if we define $\mathcal{H} = \{\|\Delta\| \leq n^{\zeta-(1/2)}\}$, then for all $n \geq n_0$, $\mathcal{H} \subseteq \mathcal{G}$. Since we assumed all moments of the principal components $\xi_j$ and the errors $\varepsilon_i$ to be finite, then Markov's inequality is readily used to show that $P(\mathcal{H}) = 1 - O(n^{-D})$ for all $D > 0$. It follows that $P(\mathcal{F}) = 1 - O(n^{-D})$, and so (5.24) implies (4.8).

5.5. *Proof of (5.14) and (5.15).* Define $\widehat{\Delta}_j$ by

$$\hat{\phi}_j(t) = \phi_j(t) + \sum_{k\,:\,k\neq j} (\theta_j - \theta_k)^{-1} \phi_k(t) \int \Delta\phi_j\phi_k + \widehat{\Delta}_j(t). \tag{5.25}$$

It may be proved that

$$\hat{\phi}_j - \phi_j = \sum_{k\,:\,k\neq j} (\hat{\theta}_j - \theta_k)^{-1} \phi_k \int \Delta\hat{\phi}_j\phi_k + \phi_j \int_{\mathcal{I}} (\hat{\phi}_j - \phi_j)\phi_j,$$

from which it follows that

$$\widehat{\Delta}_j = \sum_{k\,:\,k\neq j} \{(\hat{\theta}_j - \theta_k)^{-1} - (\theta_j - \theta_k)^{-1}\} \phi_k \int_{\mathcal{I}} \Delta\hat{\phi}_j\phi_k$$

$$+ \sum_{k\,:\,k\neq j} (\theta_j - \theta_k)^{-1} \phi_k \int_{\mathcal{I}} \Delta(\hat{\phi}_j - \phi_j)\phi_k + \phi_j \int_{\mathcal{I}} (\hat{\phi}_j - \phi_j)\phi_j.$$

If $\mathcal{F}$ holds then so too does the event $\mathcal{E}$ and, in view of (4.2), $|\theta_j - \theta_k| \leq 2|\hat{\theta}_j - \theta_k|$ for all $1 \leq j \leq m$ and all $k \neq j$. Therefore, writing $p = \sum_{j\geq 1} p_j\phi_j$ and using (5.1), we deduce that

$$\left| \int_{\mathcal{I}} p\widehat{\Delta}_j \right| \leq 2|\hat{\theta}_j - \theta_j| \left\{ \sum_{k\,:\,k\neq j} (\theta_j - \theta_k)^{-4} p_k^2 \right\}^{1/2} \|\Delta\hat{\phi}_j\|$$

$$+ \left| p_j \int_{\mathcal{I}} (\hat{\phi}_j - \phi_j)\phi_j \right| \tag{5.26}$$

$$+ \left\{ \sum_{k\,:\,k\neq j} (\theta_j - \theta_k)^{-2} p_k^2 \right\}^{1/2} \left\| \int \Delta(\hat{\phi}_j - \phi_j) \right\|.$$

Since $|p_j| \leq \text{const.} j^{-\pi}$ for each $j$ then, if $d = 2$ or $4$,

$$\sum_{k\,:\,k\neq j} (\theta_j - \theta_k)^{-d} p_k^2 \leq \text{const.} \{t_{\alpha d - 2\pi}(j) + j^{\alpha d + d - 2\pi}\}$$

$$\leq \text{const.} (1 + j^{\alpha d + d - 2\pi}).$$

Moreover, $\|\Delta\hat{\phi}_j\| \leq \|\Delta\phi_j\| + \|\Delta(\hat{\phi}_j - \phi)\|$, $E\|\Delta\phi_j\|^2 \leq \text{const.} n^{-1}\theta_j$, and if $\mathcal{F}$ holds, $\|\Delta(\hat{\phi}_j - \phi)\| \leq \text{const.} \|\Delta\|^2 \delta_j^{-1}$. We shall show in Section 5.6 that

if $\eta$, in the definition of $\mathcal{F}$ at (5.13), is chosen sufficiently small, then whenever $\mathcal{F}$ holds, $|\int_{\mathcal{I}} (\hat{\phi}_j - \phi_j)\phi_j| \leq C_0\hat{a}_j$ for $1 \leq j \leq m$, where $C_0 > 0$ is a constant depending on neither $j$ nor $n$, and $\hat{a}_j$ is a nonnegative random variable satisfying $E(\hat{a}_j^2) \leq n^{-2}j^4$. $\tag{5.27}$

Combining (5.26) and the results in this paragraph, we deduce that

$$(5.28) \quad \begin{aligned} E\bigg\{ I(\mathcal{F}) &\bigg(\int_{\mathcal{I}} p\widehat{\Delta}_j\bigg)^2 \bigg\} \\ &\leq \text{const.}\{n^{-2}j^{-\alpha}(1+n^{-1}j^{3\alpha+2})(1+j^{4\alpha+4-2\pi}) \\ &\qquad + n^{-2}j^{\alpha+1}(1+j^{2\alpha+2-2\pi}) + n^{-2}j^{4-2\pi}\}. \end{aligned}$$

Note too that

$$(5.29) \quad \begin{aligned} E\bigg\{ \sum_{k\,:\,k\neq j} &(\theta_j - \theta_k)^{-1} p_k \int \Delta\phi_j\phi_k \bigg\}^2 \\ &\leq \bigg\{ \sum_{k\,:\,k\neq j} (\theta_j - \theta_k)^{-2} p_k^2 \bigg\} E\bigg\| \int \Delta\phi_j \bigg\|^2 \\ &\leq \text{const.}\, n^{-1}j^{-\alpha}(1+j^{2\alpha+2-2\pi}). \end{aligned}$$

When $p = g$ we may substitute $\pi = \alpha + \beta$ into (5.28). Then we can deduce from (5.28) that, assuming $\alpha + 2 \leq \beta$ as well as the bound $j \leq m \leq n^{1/(\alpha+2\beta-1)}$, the right-hand side of (5.28) is bounded above by a constant multiple of $n^{-1}j^{-\alpha}$. Since $\beta > 1$ then this bound also applies to the right-hand side of (5.29).

In the case $p = x$ the fact that $\alpha + 2 \leq \beta$, as well as the bound $j \leq m \leq n^{1/(\alpha+2\beta-1)}$, imply that the right-hand side of (5.28) is dominated by the right-hand side of (5.29). Hence, for both $p = g$ and $p = x$ the bound at (5.14) follows from (5.25), (5.28) and (5.29).

Observe too that by (5.28)

$$(5.30) \quad \begin{aligned} E\bigg\{ I(\mathcal{F}) &\sum_{j=1}^m q_j \int_{\mathcal{I}} p\widehat{\Delta}_j \bigg\}^2 \\ &\leq m\,\text{const.}\{n^{-2}t_{\alpha+2\kappa+1}(m) + n^{-2}t_{3\alpha+2\kappa+4-2\pi}(m) \\ &\qquad + n^{-3}t_{2\alpha+2\kappa+2}(m) + n^{-3}t_{6\alpha+2\kappa+6-2\pi}\}. \end{aligned}$$

Now, $\kappa - \pi = -(\beta+\gamma)$ if $p = g$, and it equals $-(\alpha+\beta+\gamma)$ if $p = x$. Therefore, if $p = g$ then $3\alpha + 2\kappa + 4 - 2\pi = 3\alpha + 4 - 2(\beta+\gamma) < (\alpha+2\beta-1) - 1$, and $6\alpha + 2\kappa + 6 - 2\pi = 2\{3\alpha + 3 - (\beta+\gamma)\} < 2(\alpha+2\beta-1) - 1$. [We subtract the extra 1 to account for the factor $m$ on the right-hand side of (5.30).] These two results, and the fact that $m^{\alpha+2\beta-1} \leq n$, imply that the terms in $mn^{-2}t_{3\alpha+2\kappa+4-2\pi}(m)$ and $mn^{-3}t_{6\alpha+2\kappa+6-2\pi}$ in (5.30) may be replaced by $n^{-1}$ without affecting the validity of the bound when $p = g$. Furthermore, when $p = g$, $\alpha + 2\kappa + 1 = 3\alpha - 2\gamma + 1 < (\alpha+2\beta-1) - 1$ and $2\alpha + 2\kappa + 2 = 4\alpha - 2\gamma + 2 < 2(\alpha+2\beta-1) - 1$, and so the terms in $mn^{-2}t_{\alpha+2\kappa+1}(m)$ and $mn^{-3}t_{2\alpha+2\kappa+2}$ may also be replaced by $n^{-1}$. Therefore the right-hand of (5.30) may be replaced by $n^{-1}$ when $p = g$. An identical argument shows this also to be the case when $p = x$. Hence, in either setting,

$$(5.31) \qquad E\left\{ I(\mathcal{F}) \sum_{j=1}^{m} q_j \int_{\mathcal{I}} p\widehat{\Delta}_j \right\}^2 \leq \text{const.} n^{-1}.$$

Using (4.3) it can be proved that

$$(5.32) \qquad nE\left\{ \sum_{j=1}^{m} \sum_{k \, : \, k \neq j} (\theta_j - \theta_k)^{-1} q_j p_k \int \Delta \phi_j \phi_k \right\}^2 \leq \text{const.} t_{2\kappa - \alpha}(m).$$

Combining (5.25), (5.31) and (5.32) we obtain (5.15).

5.6. *Proof of (5.27).* It may be proved from (5.25) that $\|\hat{\phi}_j - \phi_j\|^2 = \hat{u}_j^2 + \hat{v}_j^2$, where

$$\hat{u}_j^2 = \sum_{k \, : \, k \neq j} (\hat{\theta}_j - \theta_k)^{-2} \widehat{w}_{jk}^2, \qquad \hat{v}_j^2 = \left\{ \int (\hat{\phi}_j - \phi_j) \phi_j \right\}^2$$

and $\widehat{w}_{jk} = \int \Delta \hat{\phi}_j \phi_k$. Since both $\phi_j$ and $\hat{\phi}_j$ are of unit length then $\hat{v}_j^2 = 2\{1 - (1 - \hat{u}_j^2)^{1/2}\} - \hat{u}_j^2$, which implies that

$$(5.33) \qquad \text{for all } j \geq 1, \qquad \|\hat{\phi}_j - \phi_j\|^2 \leq 2\hat{u}_j^2, \qquad \hat{v}_j^2 \leq \hat{u}_j^4.$$

If the event $\mathcal{F}$ obtains then $|\hat{\theta}_j - \theta_k|^{-1} \leq 2|\theta_j - \theta_k|^{-1}$ for all $j, k$ such that $j \neq k$ and $1 \leq j \leq m$. For the same range of values of $j$ and $k$, $|\theta_j - \theta_k|^{-1} \leq D\theta_m^{-1} m$. Here $D = C^2$, where $C$ is as in (4.2). Defining $\hat{x}_{jk} = \int \Delta \phi_j \phi_k$ and $\hat{y}_{jk} = \int \Delta(\hat{\phi}_j - \phi_j)\phi_k$, we have $\widehat{w}_{jk}^2 \leq 2(\hat{x}_{jk}^2 + \hat{y}_{jk}^2)$, and hence, assuming $\mathcal{F}$ holds, we have for $1 \leq j \leq m$,

$$(5.34) \begin{aligned} \hat{u}_j^2 &\leq 8 \sum_{k \, : \, k \neq j} (\theta_j - \theta_k)^{-2}(\hat{x}_{jk}^2 + \hat{y}_{jk}^2) \leq 8\hat{A}_j + 8D^2\theta_m^{-2}m^2\hat{c}_j \\ &\leq 8\hat{A}_j + 8D^2\theta_m^{-2}m^2\|\Delta\|^2\|\hat{\phi}_j - \phi_j\|^2, \end{aligned}$$

where $\hat{A}_j = \sum_{k \, : \, k \neq j} (\theta_j - \theta_k)^{-2}\hat{x}_{jk}^2$ and $\hat{c}_j = \sum_{k \, : \, k \neq j} \hat{y}_{jk}^2 \leq \|\Delta\|^2\|\hat{\phi}_j - \phi_j\|^2$.

Condition (4.3) implies that $nE(\hat{x}_{jk}^2) \leq \text{const.} \theta_j \theta_k$, where the constant does not depend on $j$, $k$ or $n$. Moreover,

$$\sum_{k \, : \, k \neq j} (\theta_j - \theta_k)^{-2} \theta_j \theta_k \leq \text{const.} \sum_{k \, : \, k \neq j} \left\{ \frac{\max(j, k)}{\max(\theta_j, \theta_k)|j - k|} \right\}^2 \theta_j \theta_k \leq \text{const.} j^2.$$

Therefore, $E(\hat{A}_j) \leq \text{const.} n^{-1} j^2$ for $1 \leq j \leq m$, and similar calculations show that

$$(5.35) \qquad E(\hat{A}_j^2) \leq D_1^2 n^{-2} j^4,$$

where $D_1 > 0$ depends on neither $j$ nor $n$.

Combining (5.34) with the first part of (5.33) we deduce that if $\mathcal{F}$ holds,

$$(5.36) \qquad \|\hat{\phi}_j - \phi_j\|^2 \leq 16\hat{A}_j + 16D^2\theta_m^{-2}m^2\|\Delta\|^2\|\hat{\phi}_j - \phi_j\|^2$$

for $1 \leq j \leq m$. However, if $c > 0$ is given, and if $\eta > 0$ is chosen sufficiently small in the definition of $\mathcal{F}$ at (5.13), then for all sufficiently large $m$, $\mathcal{F}$ implies $\|\Delta\| \leq cm^{-1}\theta_m$. Hence, by (5.36), if $\mathcal{F}$ holds, then for $1 \leq j \leq m$,

$$(1 - 16D^2c^2)\|\hat{\phi}_j - \phi_j\|^2 \leq 16\hat{A}_j.$$

Choosing $c$ so small that $16D^2c^2 \leq \frac{1}{2}$, we deduce that if $\mathcal{F}$ holds, then for $1 \leq j \leq m$, $\|\hat{\phi}_j - \phi_j\|^2 \leq 32\hat{A}_j$. Combining this result with (5.34), and noting the choice of $c$, we deduce that if $\mathcal{F}$ holds, then for $1 \leq j \leq m$, $\hat{u}_j^2 \leq 16\hat{A}_j$. From this property and the second part of (5.33) we conclude that if $\mathcal{F}$ holds, then for $1 \leq j \leq m$,

$$(5.37) \qquad \left|\int_{\mathcal{I}}(\hat{\phi}_j - \phi_j)\phi_j\right| \leq \hat{u}_j^2 \leq \|\hat{\phi}_j - \phi_j\|^2 \leq 32\hat{A}_j.$$

Taking $\hat{a}_j = D_1^{-1}\hat{A}_j$, where $D_1$ is as at (5.35), and letting $C_0 = 32D_1$, we see that (5.27) follows from (5.35) and (5.37).

**6. Proof of Theorem 4.2.** We shall treat only the cases $2\gamma < \alpha + 1$ and $2\gamma = \alpha + 1$, since the third setting, $2\gamma > \alpha + 1$, is relatively straightforward. For notational simplicity we shall assume that $C_1$, in the definition of $\mathcal{B}(C_1, \beta)$, satisfies $C_1 \geq 1$, and take $\theta_j = j^{-\alpha}$ and $x_j = j^{-\gamma}$. More general cases are easily addressed.

Since $X$ is Gaussian then we may write $X_i = \sum_{j \geq 1}\xi_{ij}\phi_j$ for $i \geq 1$, where the variables $\xi_{ij}$ are independent and normal with zero mean and respective variances $\theta_j$ for $j \geq 1$. Define $\nu$ to be the integer part of $n^{1/(\alpha+2\beta-1)}$, and let $B_0 \equiv 0$ and $B_1 = \sum_{\nu+1 \leq j \leq 2\nu} j^{-\beta}\phi_j$; both are functions in $\mathcal{B}(C_1, \beta)$.

Note that $T(B_0) = 0$ and that for large $n$,

$$(6.1) \qquad T(B_1) \geq \text{const.}\,n^{-(\beta+\gamma-1)/(\alpha+2\beta-1)},$$

where, here and below, "const." denotes a finite, strictly positive, generic constant. Write $\Xi_i = \sum_{\nu+1 \leq j \leq 2\nu}\xi_{ij}j^{-\beta}$. The observed data are $Y_i = t\Xi_i + \varepsilon_i$ for $1 \leq i \leq n$, where $t = 0$ or $1$ according as $b = B_0$ or $b = B_1$, respectively. Denote by $P_t$ the joint distribution of the $Y_i$'s for $t = 0$ or $1$. Elementary calculations show that the chi-squared distance between $P_0$ and $P_1$ is given by

$$d(P_0, P_1) = \int \frac{(dP_1)^2}{dP_0} = \exp\left(\sigma^{-2}\sum_{i=1}^n \Xi_i^2\right),$$

where $\sigma^2$ denotes the variance of the error distribution.

The variables $\Xi_i$ are independent and normally distributed with zero means and variance $V_n$, where $nV_n = n\sum_{\nu+1 \leq j \leq 2\nu} j^{-\alpha-2\beta} \to$ const. as $n \to \infty$. Indeed,

$$(6.2) \qquad\qquad E_1\{d(P_0, P_1)\} \to \text{const.,}$$

where $E_t$ denotes expectation in the model with $b = B_t$, for $t = 0$ or $1$. Let $\widehat{T}$ be any estimator such that for some $D > 0$,

$$(6.3) \qquad\qquad E_0\{\widehat{T} - T(B_0)\}^2 \leq Dn^{-2(\beta+\gamma-1)/(\alpha+2\beta-1)}.$$

Put

$$\rho = \frac{2[E_0\{\widehat{T} - T(B_0)\}^2 E_1\{d(P_0, P_1)\}]^{1/2}}{|T(B_1) - T(B_0)|}.$$

It follows from (6.1), (6.2) and the fact that $T(B_0) = 0$, that if $D$ in (6.3) is chosen sufficiently small, $\rho \leq \frac{1}{2}$. In this case,

$$
\begin{aligned}
(6.4) \qquad E_1\{\widehat{T} - T(B_1)\}^2 &\geq \{T(B_1) - T(B_0)\}^2(1 - \rho) \\
&\geq \text{const.}\, n^{-2(\beta+\gamma-1)/(\alpha+2\beta-1)},
\end{aligned}
$$

where the first inequality follows from the constrained-risk lower bound of Brown and Low [4], and the second uses (6.1) and the property $T(B_0) = 0$. Consequently, writing $E_b$ for expectation when the slope function is $b \in \mathcal{B}$, for any estimator $\widehat{T}$

$$\sup_{b \in \mathcal{B}} E_b\{\widehat{T} - T(b)\}^2 \geq \max_{t=0,1} E_t\{\widehat{T} - T(B_t)\}^2 \geq \text{const.}\, n^{-2(\beta+\gamma-1)/(\alpha+2\beta-1)}.$$

The case $2\gamma = \alpha + 1$ may be treated similarly, by taking $\nu = (n/\log n)^{1/(\alpha+2\beta-1)}$ and replacing $n$ by $n/\log n$ in (6.1), (6.3) and (6.4).

**Acknowledgment.** This work was done while Tony Cai was visiting the Mathematical Sciences Institute of the Australian National University.

Department of Statistics
The Wharton School
University of Pennsylvania
Philadelphia, Pennsylvania 19104-6340
USA
E-mail: tcai@wharton.upenn.edu

Centre for Mathematics
  and Its Applications
Australian National University
Canberra, ACT 0200
Australia
E-mail: Peter.Hall@maths.anu.edu.au